\documentclass[]{article}
\usepackage[utf8]{inputenc}
\usepackage{amsfonts}
\usepackage{nccbbb}
\usepackage{graphicx}
\usepackage{multirow}

\def\Zp{\bbbz_{+}}

\begin{document}
\title{Numerical experiments in the problems of asymptotic representation theory}
\author{Anatoly Vershik \and Dmitry Pavlov}

\maketitle

\begin{abstract}
The article presents the results of experiments in computation of
statistical values related to Young diagrams, including the estimates
on maximum and average (by Plancherel distribution) dimension of
irreducible representation of symmetric group $S_n$. The computed
limit shapes of two-dimensional and three-dimensional diagrams
distributed by Richardson statistics are presented as well.
\end{abstract}

\section{Introduction}
A classical definition of Young diagram\cite{fulton}
is the following: the Young diagram if size $n$ is a
finite descending ideal in the lattice $\Zp\times\Zp$,
i.e. a set of cells in positive quadrant, which contains
a cell $(i,j)$ only if it contains as well all the
sells under $(i,j)$ in terms of particular order.

Standard Young table of size $n$ is a Young diagram whose
cells are filled with numbers in the range from $1$
to $n$, which increase in each row and each column.
In other words, Young table is a path in a lattice of
Young diagrams, which starts at the empty (zero-sized) diagram
and ends in the given diagram.

Young diagrams correspond to irreducible representations
of $S_n$ (see for example~\cite{fulton,vershik-okounkov}),
and the number of Young tables that fit into a given Young
diagram is equal to the dimension of the corresponding
irreducible representation. For the sake of brevity we call
this number the dimension of the diagram $\Lambda$ and
denote it as $\dim(\Lambda)$.

This paper is organized as follows: in section~2 we will
present the results of numerical experiments on asymptotics
of a typical (by Plancherel measure) dimension of irreducible
representation of the symmetric group. This measure
was introduced in~\cite{vershik-kerov-77}, for further
details see~\cite{vershik-kerov}.

Section~3 is devoted to the similar computations of
a maximum dimension of irreducible representation of
the symmetric group.

The results of section~2 should be considered as a support of
a conjecture~\cite{vershik-kerov} about the existence of a
limit value of normalized dimension of a typical $\Lambda$
by Plancherel measure. This hypothetical
limit value is denoted as $c$ and called in ~\cite{vershik-kerov}
``specific entropy of an irreducible representation''.

At the same time, basing on out computations in section~3,
nothing conclusive can be said about the asymptotic behavior
of the maximum dimension. That is, there is no claim that
$c_n$ comes to plateau on approachable values of $n$.

In section~4, we are concerned with another distribution
on Young diagrams: Richardson distribution, which was studied
in~\cite{rost}. We give an experimental evidence of
the proved theorem about a limit shape of typical
Young diagrams. Then we make the similar experiments
in three dimensions (this time the dimension has nothing to do
with $S_n$), and form a hypothesis on a 3-dimensional limit
shape.

\section{Asymptotic behavior of a typical dimension of
irreducible representation of $S_n$ by Plancherel measure}
Let $\hat{S_n}$ be the collection of equivalence classes
of complex irreducible representations of $S_n$.
For $\Lambda_n\in\hat{S_n}$, we denote by $\dim\Lambda_n$
the dimension of the representation $\Lambda_n$, and by
\begin{equation}\label{mu-lambda}
\mu_n(\Lambda_n) = \frac{\dim^2 \Lambda_n}{n!}
\end{equation}
we denote the Plancherel measure~\cite{vershik-kerov-77}.
This is actually a probabilistic measure on $\hat{S_n}$,
which we can derive from Burnside's conjecture:
$$\sum_{\Lambda_n\in\hat{S_n}} \dim^2\Lambda_n = n!$$

The set $\hat{S_n}$ and the dimension $\dim\Lambda_n$ have
the following interpretation~\cite{fulton}:
$\hat{S_n}$ is the ensemble of all Young diagrams of size $n$,
and $\dim\Lambda_n$ is the number of standard Young tables that
fit into the diagram $\Lambda_n$. In the present article,
we make no difference between  $\Lambda_n$ and the corresponding
irreducible representation.

The following normalization of a dimension of diagrams with $n$
cells (see~\cite{vershik-kerov}) is a key to studying the
asymptotics of the dimension with $n\to\infty$:

\begin{equation}\label{c-lambda}
c(\Lambda_n) = {-2\over\root\of n}\ln{\dim\Lambda_n\over\root\of{n!}}
\end{equation}
We will call the coefficient $c(\Lambda_n)$ the \textsl{normalized dimension}.
From~\cite{vershik-kerov}, the following two-sided estimates are known for
$c(\Lambda_n)$:

$$\lim_{n\to\infty} \mu_n\left\{\Lambda_n:
   c_0< c(\Lambda_n) < c_1 \right\} = 1.$$
$$(c_0=\frac{2}{\pi}-\frac{4}{\pi^2}\approx0.2313,
  \ c_1=\frac{2\pi}{\sqrt 6}\approx2.5651).$$

In other words, asymptotically almost all diagrams have the dimension
lying inside the following range:
$$\root\of {n!} e^{-{c_1\over 2}\root\of n} < \dim\Lambda_n <
  \root\of {n!} e^{-{c_0\over 2}\root\of n}.$$

Vershik and Kerov~\cite{vershik-kerov} stated the following conjecture:
there exists a limit $\lim_{n\to\infty}c(\Lambda_n)$ almost everywhere
(by Plancherel measure) on the set of infinite Young tables. Infinite
Young table is an infinite sequence of nested Young  diagrams of
sizes increasing one-by-one:
$\Lambda_1 \subset \Lambda_2 \subset \Lambda_3 \subset \dots$.
The hypothetical limit value is denoted as $c$ and called
``specific entropy of an irreducible representation''.

We are studying the behavior of coefficient $c(\Lambda_n)$ w.r.t. $n$.
Some experiments in this area were undertaken
in~\cite{vershik-gribov-kerov}, where the expectation and variance
of $c(\Lambda_n)$ are given for five values of $n$, 1600 maximum
(for $n=1600$ the sample was only 14 diagrams). We have done
the same computations on modern hardware, which is far more powerful
than  25 years ago.

\bigskip
\noindent\textbf{Note.}
After we did the described experiments, it was reported of a new
paper by Alexander Bufetov\footnote{Alexander I. Bufetov.
On the Vershik-Kerov Conjecture Concerning the Shannon-Macmillan-Breiman
Theorem for the Plancherel Family of Measures on the Space of
Young Diagrams. arXiv:1001.4275 },
which gives a proof
of existence of a limit $C$ of $c(\Lambda_n)$ in $L^2(Y)$
by Plancherel measure:
$$\lim_{n\to\infty}\int(c(\Lambda_n)-C)^2 d(\mu_n(\Lambda_n)) = 0$$

The actual value of $C$ is not likely to be obtainable by methods used
in Bufetov's proof. In should be also noted that the conjecture~\cite{vershik-kerov}
about the existence of a limit if $c(\Lambda_n)$
almost everywhere is still open, and if it happens that such a limit
value $c$ exists, then clearly $C=c$.

The next two sections of the paper will describe well-known helper
routines: generating the random diagrams by Plancherel measure with
RSK algorithm and counting the dimension of the diagram.

\subsection{Hook formula}
Hook formula (see for example ~\cite{fulton,knuth3}) allows to compute
$\dim\Lambda$, avoiding the enumeration of all Young tables that fit
into $\Lambda$.
\begin{equation}\label{hook-formula}
\dim\Lambda = \frac{n!}{\prod_{(i,j)\in \Lambda}h_Y(i,j)},
\end{equation}
where $(i,j)$ is a cell of Young diagram, and $h_Y(i, j)$ is the length
of the hook associated to that cell. The hook of the cell $(i,j)$
consists of the cell itself and all cells that are in $j$-th row
to the right from $(i,j)$ or in $i$-th column above $(i,j)$ (fig.~\ref{hooks}).

\begin{figure}[ht]
\begin{center}
\begin{tabular}{|c|cccc}
\cline{1-1}
1 & & & & \\
\cline{1-4}
5 & \multicolumn{1}{|c|}{3} & \multicolumn{1}{|c|}{2} & \multicolumn{1}{|c|}{1} & \\
\cline{1-5}
7 & \multicolumn{1}{|c|}{5} & \multicolumn{1}{|c|}{4} & \multicolumn{1}{|c|}{3} &
\multicolumn{1}{|c|}{1}\\
\hline
\end{tabular}
\end{center}
\caption{Young diagram hook lengths}
\label{hooks}
\end{figure}

\subsection{Generating RSK-random diagrams}
While we have the hook formula in hand, the straight-forward
generation of random diagrams according
to Plancherel distribution would be
very computationally expensive. The Robinson-Schensted-Knuth
correspondence (RSK) and the row-insertion algorithm are of big help here.

The RSK algorithm~\cite{fulton} takes as input an arbitrary
permutation  $s\in S_n$, performs a sequence of row-insertions,
and produces a pair of a standard
Young tables $(P, Q)$ whose diagrams are equal; furthermore,
there is one-to-one correspondence (RSK) between such pairs and
permutations. Hence, the uniform distribution (Haar measure)
on $S_n$ transforms into Plancherel measure on the set of
left (or right) Young tableaux. Applying the RSK correspondence
to a random permutation in $S_n$, and taking the Young diagram
$Y(P)$ of the left Young tableau in the pair, we will get
a random Young tableau according to Plancherel measure.

\subsection{Results}
For $n\leq120$ the expectancy of $c(\Lambda_n)$ for Plancherel measure
were computed precisely by the formula,
$$c_n = \sum_{\Lambda_n\in\hat{S_n}}c(\Lambda_n)\frac{\dim^2\Lambda_n}{n!}.$$
The results are listed in the next section: see table~\ref{c-max-avg} and
figure~\ref{c-max-avg-plot}.

\begin{figure}[ht]\begin{center}
\includegraphics{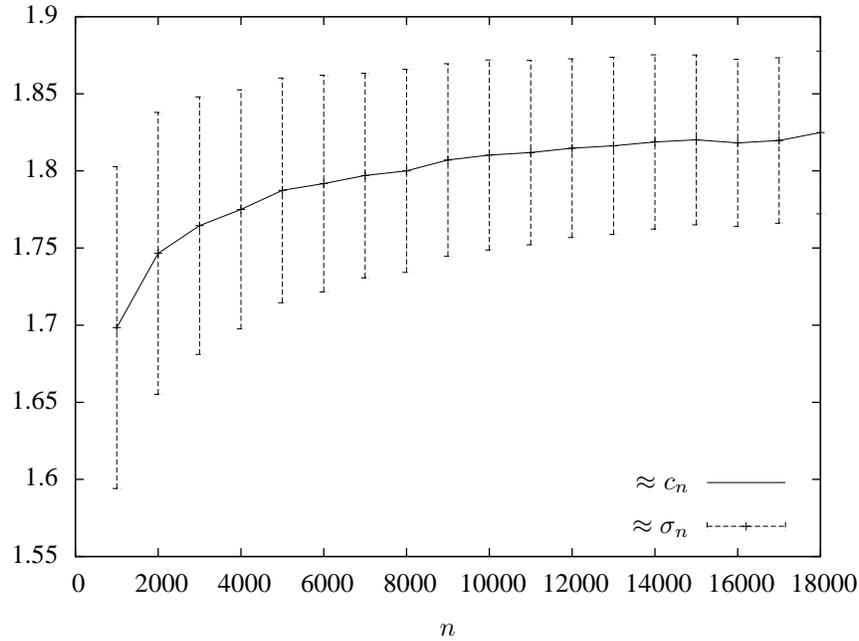}
\caption{Expected values of $c(\Lambda_n)$.
At each $n$, the height of a vertical line is equal to $\sigma_n$.}
\label{plancherel-c-avg}
\end{center}\end{figure}

There are 1,844,349,560 Young diagrams of size 120. For larger $n$,
the expectancies of $c(\Lambda_n)$ were calculated by Monte-Carlo
method, using a sample of RSK-random diagrams.
The normalized dimension of each diagram was computed by
formula~(\ref{c-lambda}), using the hook formula~\ref{hook-formula}
for obtaining $\dim\Lambda$. The procedure was run against various
$n$ in the range from 1000 to 18000. The expectation value $c_n=E(c(\Lambda_n))$
and the standard deviation $\sigma_n=\sigma(c(\Lambda_n))$ of $c(\Lambda_n)$
are listed in table~\ref{plancherel-c} and figure~\ref{plancherel-c-avg}.
From these results we can suggest that $c_n$ is asymptotically
increasing and allegedly has a limit. To be quite honest, we should notice
that the $c_n$ does not monotonically increase in the selected range: for example,
it decreases from $n=15000$ to $n=16000$. This fact was re-checked and approved
with a sample with $40000$ items. The 3-rd decimal place remained
constant after the amount of items had reached $20000$.

\begin{table}[ht]
\begin{center}
\begin{tabular}{|c|c|c|c|}
\hline
$n$ & sample size & $\approx c_n$ & $\approx\sigma_n$ \\
\hline
1000   & 2000   & 1.6984314     &  0.10431497  \\
2000   & 2000   & 1.746588      &  0.091339454 \\
3000   & 2000   & 1.7644972     &  0.08351989  \\
4000   & 2000   & 1.7750576     &  0.07747431  \\
5000   & 2000   & 1.7873781     &  0.07282907  \\
6000   & 2000   & 1.7917556     &  0.07022077  \\
7000   & 2000   & 1.7969893     &  0.06630529  \\
8000   & 2000   & 1.8000197     &  0.06586118  \\
9000   & 2000   & 1.8070668     &  0.06243244  \\
10000  & 10000  & 1.8102994     &  0.061589677 \\
11000  & 10000  & 1.8118591     &  0.059796795 \\
12000  & 10000  & 1.8147597     &  0.057941828 \\
13000  & 10000  & 1.8162445     &  0.05743194  \\
14000  & 10000  & 1.8187699     &  0.056453623 \\
15000  & 20000  & 1.820125      &  0.05504108  \\
16000  & 20000  & 1.8181555     &  0.054255717 \\
17000  & 20000  & 1.8197316     &  0.053651392 \\
18000  & 20000  & 1.8249108     &  0.052745327 \\
\hline
\end{tabular}
\end{center}
\caption{Expected values and standard deviation of $c(\Lambda_n)$}
\label{plancherel-c}
\end{table}

\subsection{Individual evolution of dimension of typical diagram}
We define the Plancherel measure on infinite Young tableaux as
a Markov measure having the following property: the corresponding
measure on Young diagrams for each $n$ is Plancherel measure
(\ref{mu-lambda}) on diagrams. It is not a big problem to derive
the transition probabilities for this Markov measure. Given that
tables which have equal diagrams have as well equal measures,
we can say that the measure of $\lambda$-shaped Young tableau is
$\frac{\dim\lambda}{n!}$. Therefore, the probability of transition
from $\lambda$ to $\Lambda$ is
\begin{equation}\label{transition}
P(\Lambda|\lambda)=\frac{\dim(\Lambda)}{(n+1)\dim(\lambda)}.
\end{equation}
(See for example~\cite{vershik-tsilevich}). So, the
conjecture~\cite{vershik-kerov} on existence of
limit of $c(\Lambda_n)$ almost everywhere means that
for almost every infinite Young tableaux
$\{\Lambda_n, n=1,2,\dots\}$, generated by the described Markov
process, $c(\Lambda_n)$ will converge to some common limit value, which obviously
will be equal to the limit value for expectation of normalized
dimension. Using the formula~(\ref{transition}) for transition probability,
we simulate the Markov process and obtain the sequence of Young diagrams
of increasing size, each according to Plancherel measure, and each
containing all the previous..

Our experiments have shown very chaotic behavior of the
normalized dimension of such sequences, which probably
points out that the Vershik-Kerov conjecture on
existence of limit of $c(\lambda_n)$ almost everywhere
is not easy to prove nor to prove wrong.

\begin{figure}[ht]\begin{center}
\includegraphics[scale=0.8]{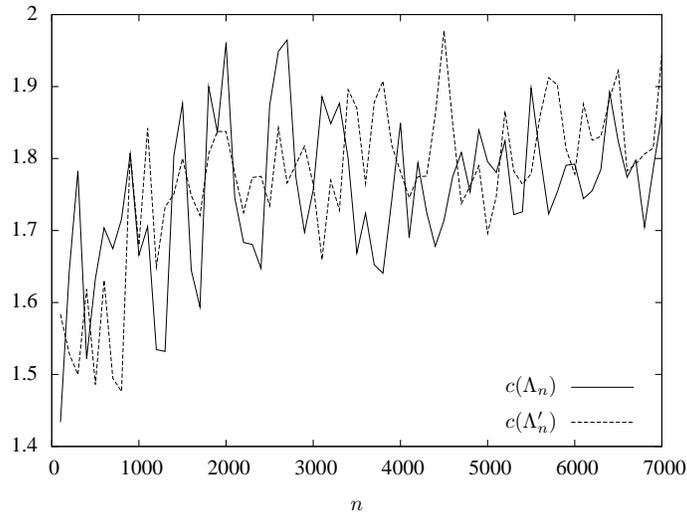}
\caption{$c(\lambda_n)$ of two random sequences of Young diagrams
by Plancherel measure}
\label{single-plancherel}
\end{center}\end{figure}

Non-regularity of behavior of normalized dimension
is illustrated on figure ~\ref{single-plancherel},
which depicts the values of
$c(\lambda_n)$ for two Markov sequences of Young diagrams.
The values are computed for $n\in[100..7000]$ (only multiples
of 100 were taken).

\section{Asymptotic behavior of a maximum dimension of
irreducible representation of $S_n$}
In this section we will study the behavior of {\it maximum dimension}
of a diagram of size $n$
$$m_n=\max_{\Lambda_n\in\hat{S_n}}\dim\Lambda_n$$
\noindent and its normalized value
$$\overline{c_n}=c(\overline\Lambda_n),$$
where $\overline\Lambda_n$ is the diagram of size $n$, which has
a maximum dimension over all diagrams of size $n$.

The problem of computing maximum dimension was stated
in 1968 (see~\cite{baer}). In a paper by McKay~\cite{mckay}
there are values of $\max\dim\Lambda_n$ for $n$ up to 75.
Basing on his results, McKay assumed that
\begin{equation}\label{mckay-assumption}
\max\frac{\dim\lambda_n}{\sqrt{n!}}\leq\frac{1}{n}.
\end{equation}
This assumption was the opposite to an alternative
hypothesis, stating that there are irreducible representations
of arbitrary large dimension, for which the
inequality~\ref{mckay-assumption} is not true.
Right before his paper was published, McKay sadly admitted that
the alternative hypothesis is true for $n=81$. Nevertheless,
as shown in~\cite{vershik-kerov}, McKay's assumption is asymptotically
true, and even stronger fact is true: $\frac{\max\dim\Lambda_n}{\sqrt{n!}}$
with  $n\to\infty$ is decreasing as $e^{-c\sqrt n}$, i.e.~not just faster
than $1/n$, but faster than any polynomial fraction.
The estimates on normalized dimension, given in~\cite{vershik-kerov}
for typical Young diagram, are the same for maximum dimension,
and while both have the save logarithmic order, the experiment
shows that the constants are different.

For $n$ up to $130$ we find the $\max\dim\Lambda_n$ via enumeration
of all Young diagrams of size $n$. In fact, we enumerated not
diagrams, but \textit{partitions} of $n$, using then the trivial
correspondence between Young diagrams and partitions of
integer~\cite{fulton}. The dimension of each diagram
was computed using hook formula~\ref{hook-formula}.

There are 5,371,315,400 Young diagrams of size $n=130$.
For larger $n$, because of computational inability to
enumerate all Young diagrams, the set of diagrams was
restricted: we considered only symmetric diagrams,
\textit{or diagrams that can be obtained from symmetric
by adding one cell}. This restriction often does not affect
the final result, but for example for $n=14$ the diagram
with maximum dimension does not pass.
However, this restriction makes no big noise after all.
Table~\ref{c-max-dim} contains values of $\overline{c_n}$.
For $n=310$, 151,982,627 diagrams were enumerated.

\begin{table}[htb]\begin{center}
\begin{tabular}{|c|c|c|c|}
\hline
$n$ & $\overline{c_n}$ & $n$ & $\approx\overline{c_n}$ \\
\hline
10      & 0.57453286  & 140  & 1.05010306 \\
20      & 0.8198125   & 150  & 1.0839802  \\
30      & 0.7912792   & 160  & 1.05304872 \\
40      & 0.86301332  & 170  & 1.0784368  \\
50      & 0.90097636  & 180  & 1.0775954  \\
60      & 0.94780416  & 190  & 1.0940416  \\
70      & 0.98343194  & 200  & 1.0953336  \\
80      & 0.96466594  & 210  & 1.1026434  \\
90      & 0.9749938   & 220  & 1.11596048 \\
100     & 1.035376    & 230  & 1.1106038  \\
110     & 1.02168428  & 240  & 1.1273114  \\
120     & 1.02246392  & 250  & 1.11251032 \\
130     & 1.0514124   & 260  & 1.11878812 \\
        &             & 270  & 1.1175388  \\
        &             & 280  & 1.1173389  \\
        &             & 290  & 1.13589692 \\
        &             & 300  & 1.12641788 \\
        &             & 310  & 1.148327   \\
\hline
\end{tabular}
\end{center}
\caption{Values of $c(\Lambda_n)$ for maximum-dimension diagrams.
The first column has exact values, while the second column
has the values obtained with enumeration of restricted set of diagrams.}
\label{c-max-dim}
\end{table}

\begin{figure}[htb]\begin{center}
\includegraphics[scale=0.8]{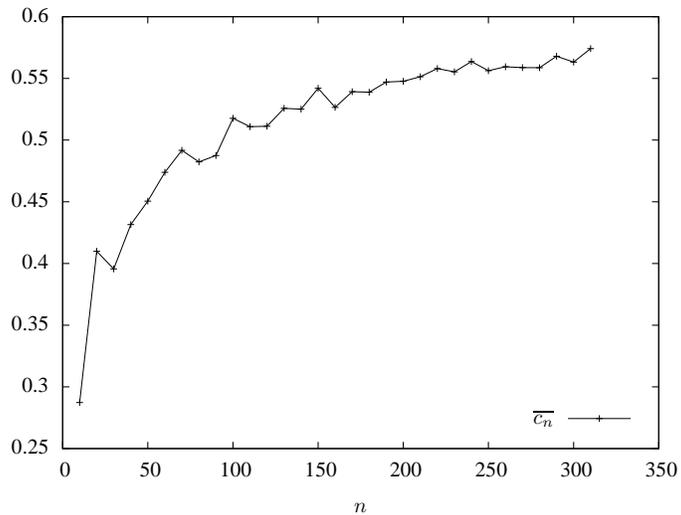}
\caption{Values of $c(\Lambda_n)$ for maximum-dimension diagrams.
The values starting from $n=140$ are approximate, because
of the mentioned restriction.}
\label{c-max-dim-plot}
\end{center}\end{figure}

While the two-sided estimates, given by Vershik and Kerov~\cite{vershik-kerov}
for maximum dimension, are the same as for typical dimension, these two values
have different behavior, and the limit of the sequence $\overline{c_n}$
is not likely to exist (see figure~\ref{c-max-dim-plot}). We also emphasize
that the the maximum dimension is way greater than the typical one
(vice versa for normalized values, because of the minus sign
in the exponent).

\begin{figure}\begin{center}
\includegraphics{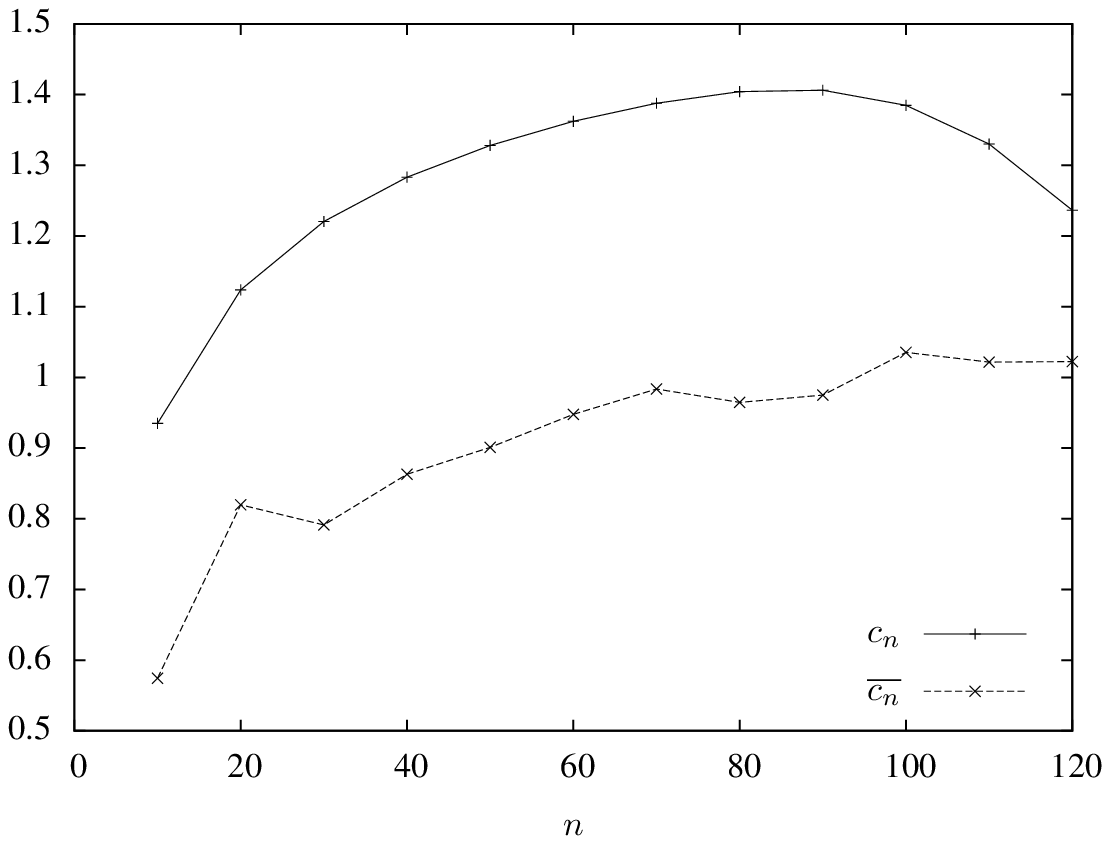}
\caption{Normalized dimensions: maximum and typical.}
\label{c-max-avg-plot}
\end{center}\end{figure}

Table~\ref{c-max-avg} and figure~\ref{c-max-avg-plot}
present the comparison of exact values of
$\overline{c_n}$ and $c_n$. Neither of functions
is monotonous.

\begin{table}
\begin{center}
\begin{tabular}{|c|c|c|}
\hline
$n$ & $\overline{c_n}$ & $c_n$ \\
\hline
10     & 0.57453287   &   0.9348365 \\
20     & 0.81981254   &   1.1238908 \\
30     & 0.7912792    &   1.2205664 \\
40     & 0.8630133    &   1.283057  \\
50     & 0.90097636   &   1.3281072 \\
60     & 0.94780415   &   1.3622344 \\
70     & 0.98343194   &   1.3878295 \\
80     & 0.96466595   &   1.4042087 \\
90     & 0.9749938    &   1.4061089 \\
100    & 1.035376     &   1.3848866 \\
110    & 1.0216843    &   1.3299882 \\
120    & 1.0224639    &   1.2363929 \\
\hline
\end{tabular}
\end{center}
\caption{Normalized dimensions: maximum and typical.}
\label{c-max-avg}
\end{table}

\section{Random diagrams by Richardson}
Rost~\cite{rost} considers a Markov process of a particle
in $\{0, 1\}^\bbbz$, which can be treated as a process
of increasing the Young diagram cell-by-cell, starting
from an empty diagram. The transition to the next state
(increasing the diagram by one cell) is performed in the
following way:

Among all diagrams of size $n+1$, which contain the
given diagram, one diagram is picked up randomly;
each in the list with equal probability.
In other words, from all the $k$ ``dimples'' of
the diagram of size $n$ one dimple is picked with
probability $1/k$. This growth process was introduced
by Richardson~\cite{richardson}. Rost~\cite{rost} found
and proved the limit shape for Young diagram w.r.t.~Richardson
measure (see below).

\begin{figure}[htb]\begin{center}
\includegraphics[scale=0.8]{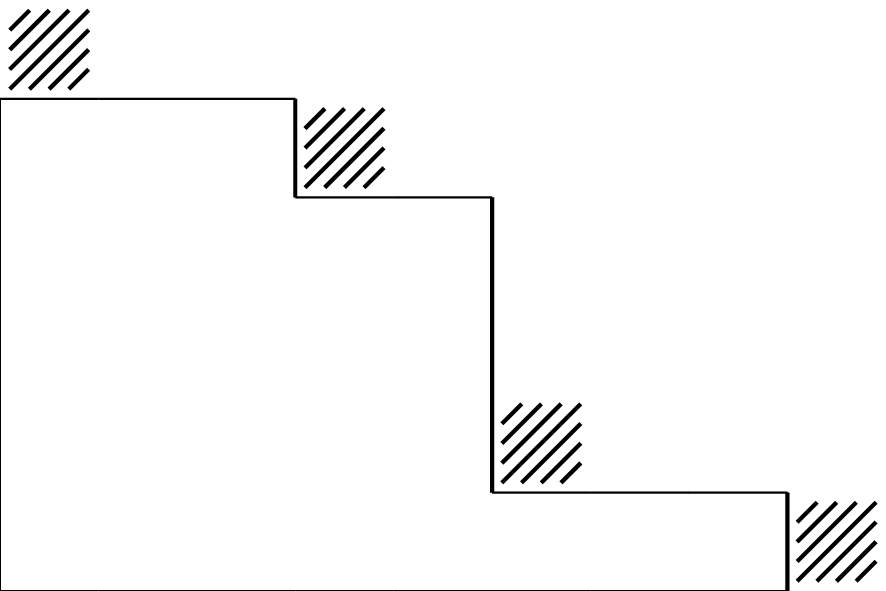}
\caption{Building of a random Young diagram by Richardson measure.
From $k$ dashed dimples, each has the probability of $1/k$.}
\end{center}\end{figure}

\begin{figure}[htb]\begin{center}
\includegraphics[scale=0.7]{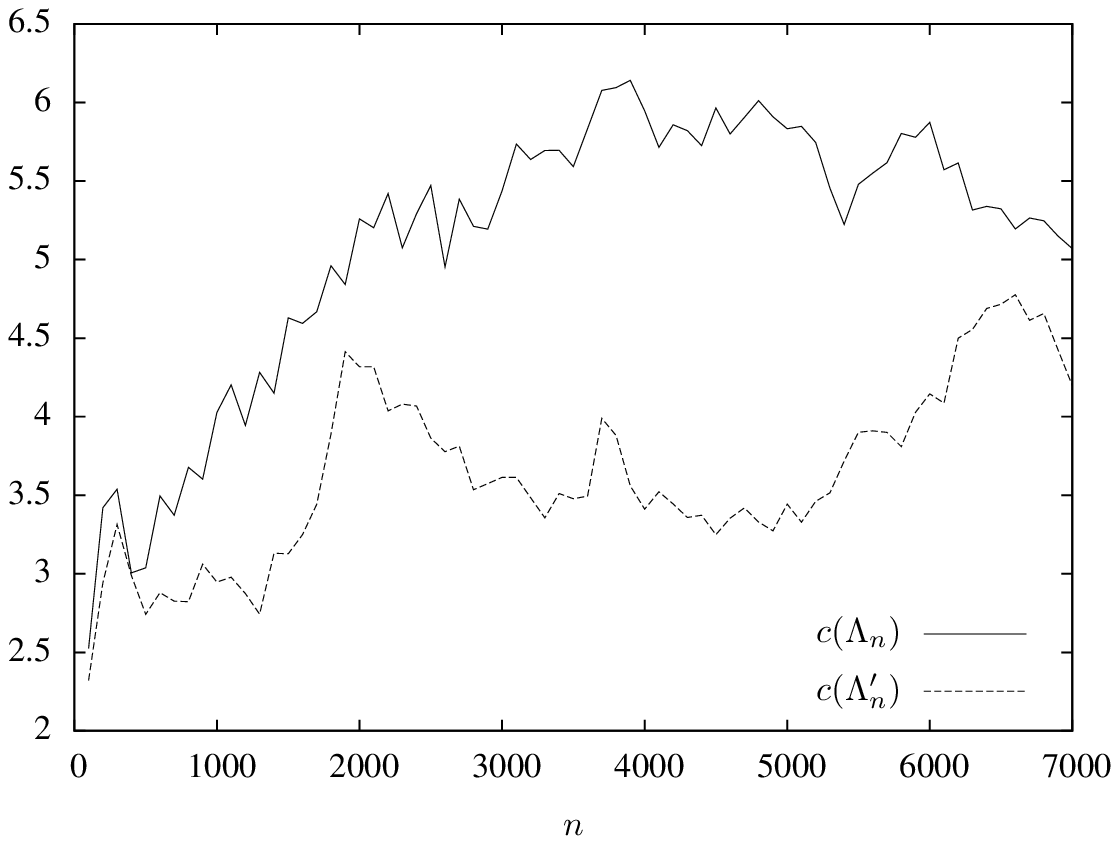}
\caption{$c(\lambda_n)$ of two random sequences of Young diagrams
by Richardson measure.}
\label{single-richardson}
\end{center}\end{figure}

We computed the values of normalized dimension
$c(\Lambda_n)$ for sequences of nested Young diagrams,
generated by Richardson process (see figure~\ref{single-richardson}).
The difference between these values and the values for Plancherel measure
(figure~\ref{single-plancherel}) makes it clear that these measures
are totally different.

Let us consider a process of infinite growth of Young diagram,
with normalization on each step to along each axis;
the area of the diagram thus remains constant.
As shown in~\cite{rost}, with $n\to\infty$ the process converges
to a limit shape, defined by equation $\sqrt{x}+\sqrt{y}=h$.
The exact value of $h$ is up to normalization.
In~\cite{rost} $h$ is equal to one: $\sqrt{x}+\sqrt{y}=1$.
Another normalization that makes sense is the normalization
by the area of the resulting figure, which is equal to $1/6$
of the area of circumscribed square. The side of the square is
equal to $h^2$, so
$$S = \int_0^{h^2}(h-\sqrt{x})^2 dx = h^4/6$$
If we take the area $S$ for 1, then the value of $h$ is
equal to $\root4\of 6$, and the limit-shape equation will be

\begin{equation}\label{sqrt-plus-sqrt}
\sqrt{x}+\sqrt{y}={\root 4 \of 6}
\end{equation}

\subsection{$d$-dimensional Young diagrams}
We call a $d$-dimensional Young diagram a finite descending ideal
in the lattice $(\Zp)^d$. Unless specified otherwise,
just ``Young diagram'' will mean a two-dimensional Young diagram.

Vershik and Kerov~\cite{vershik-kerov-77} introduced
a convenient coordinate system for presentation of
Young diagrams: the diagram is rotated by $45^\circ$
from the so-called French depiction, which conforms
to Descartes coordinates.

Similarly, $d$-dimensional Young diagrams can be represented
as \textit{functions}, defined on a $(d-1)$-dimensional
hyperplane, crossing the origin and orthogonal to the
main diagonal. The value of the function is the
length of a line span, parallel to the main diagonal,
starting at the hyperplane and ending at the border of
the Young diagram.

Having any Young diagram defined by this function,
we easily define the \textit{average shape} of
a collection of diagrams, as the average of corresponding functions.
This definition trivially applies to multi-dimensional Young
diagrams as well.

Despite the function is defined in ``rotated'' coordinate system,
we still depict the Young diagrams and their average shapes
in Descartes coordinates, by applying an inverse transformation.

\begin{figure}[htb]\begin{center}
\includegraphics[scale=0.65]{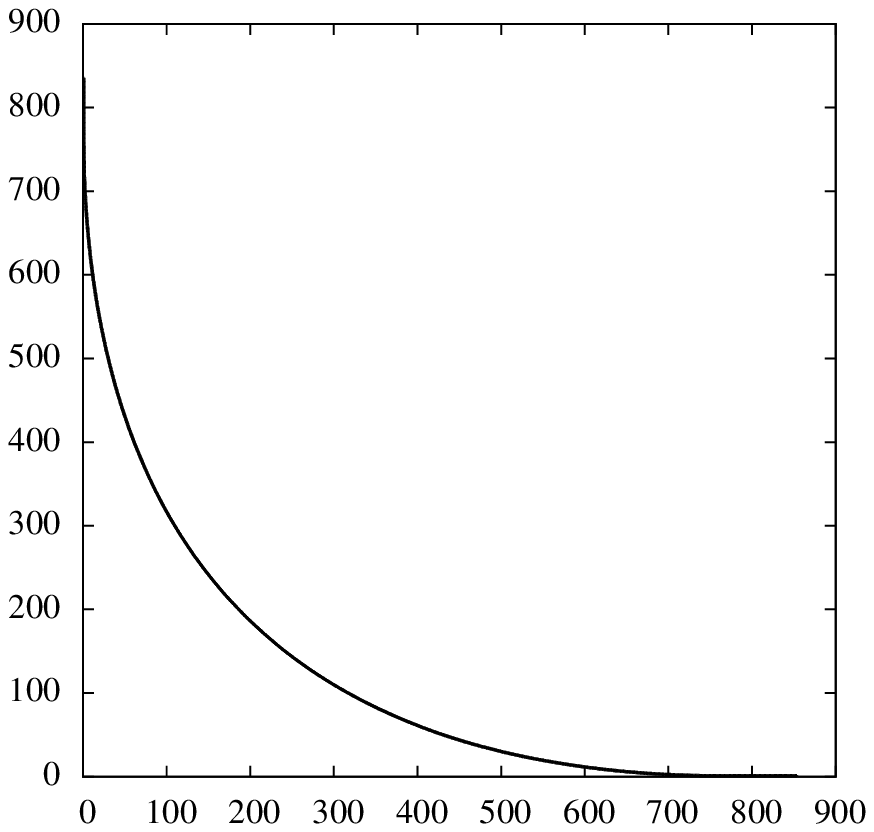}
\caption{Average shape of 2200 diagrams of size 100000}
\label{shape}
\end{center}\end{figure}

\begin{figure}[htb]\begin{center}
\includegraphics[scale=0.65]{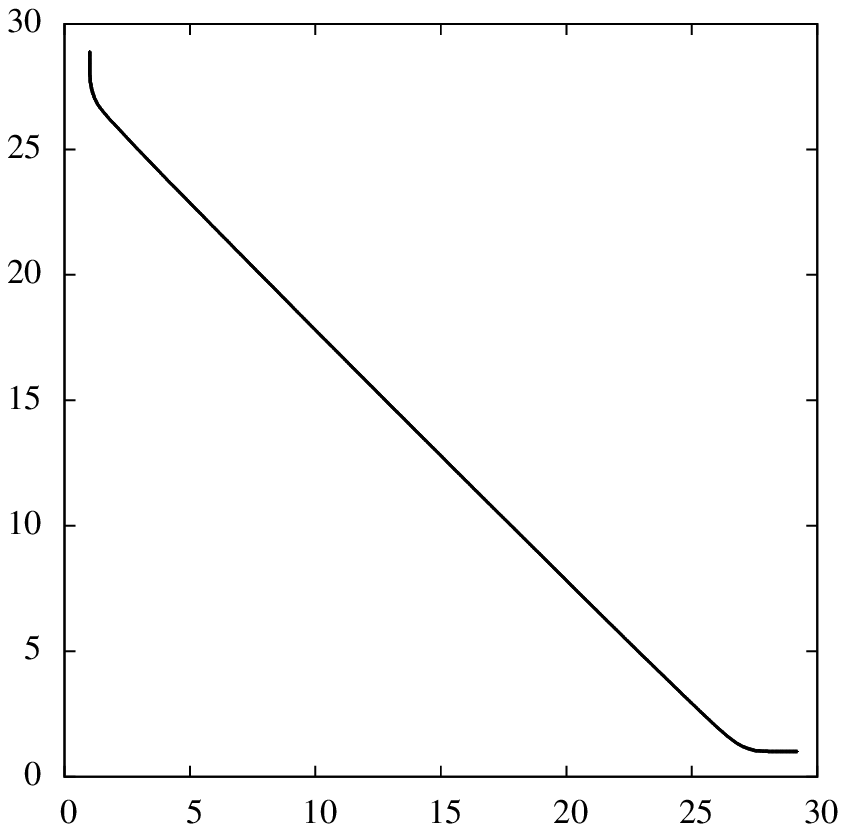}
\caption{Average shape of 2200 diagrams of size 100000
in coordinates $(\sqrt{x}, \sqrt{y})$}
\label{shape-sqrt}
\end{center}\end{figure}

The average shape of 2200 random diagrams of size $n=100000$
is shown on figure~\ref{shape}. A visual verification
of Rost's theorem can be shown by plotting the average
shape in coordinates $(\sqrt{x},\sqrt{y})$ (see figure~\ref{shape-sqrt}).

\begin{figure}[htb]\begin{center}
\includegraphics[scale=0.65]{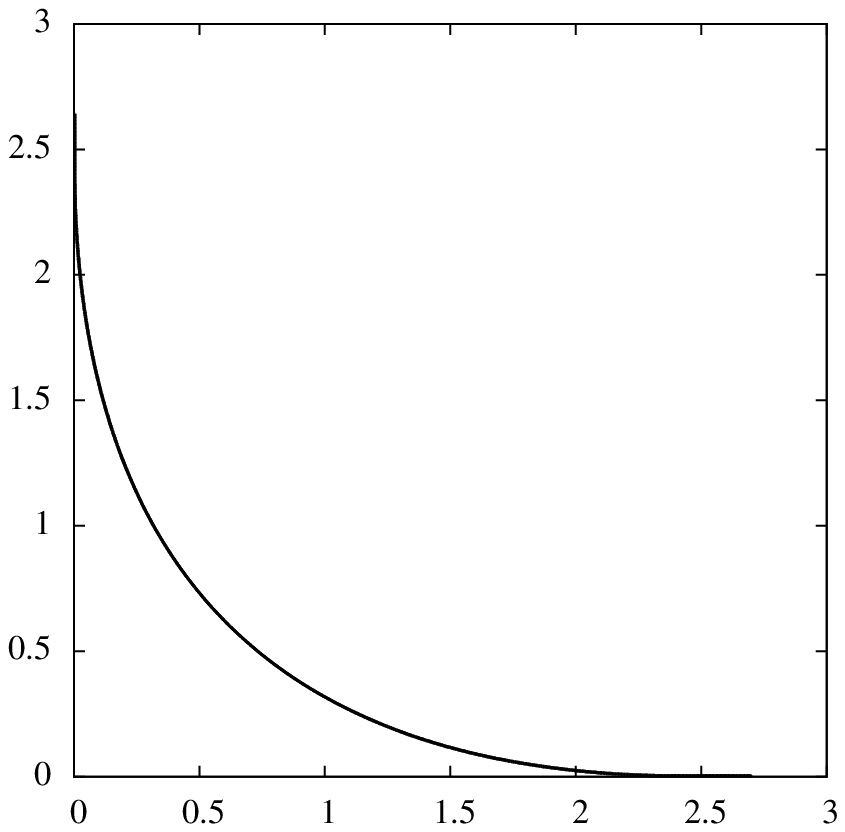}
\caption{Average shape of 2200 diagrams of size 100000,
scaled by $1/\sqrt{n}$}
\label{shape-norm}
\end{center}\end{figure}

Scaling this shape by $1/\sqrt{n}$, we get the close approximation
of the plot of equation~\ref{sqrt-plus-sqrt}. The area of
the figure~\ref{shape-norm} is equal to 1.

\subsection{Standard deviation of main diagonal segment}
\begin{figure}[htb]\begin{center}
\includegraphics[scale=0.65]{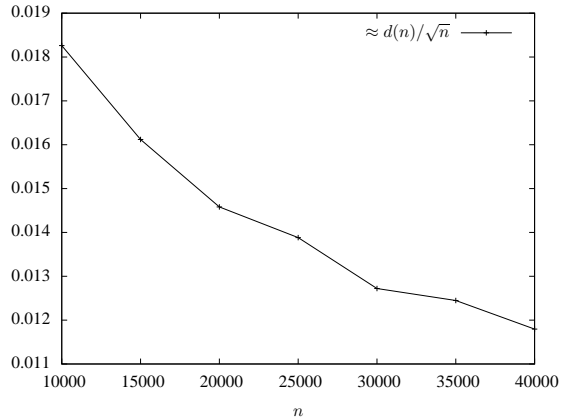}
\caption{Normalized standard deviation of main diagonal segment}
\label{deviation-norm}
\end{center}\end{figure}

In the previous section we verified that the average
shape comply with equation~\ref{sqrt-plus-sqrt}.
For verifying that the average shape is indeed a limit shape,
we computed the standard deviation of s-called
\textit{main diagonal segment}---the length of a line span
of the main diagonal, starting in origin and ending in the
average shape. In table~\ref{deviation-center} the values
of the standard deviation $d(n)$ are listed, for $n$
from 10000 to 40000, along with the normalized values
$d(n)/\sqrt{n}$. Figure~\ref{deviation-norm} shows the
decrease of the normalized standard deviation with $n\to\infty$.

\begin{table}[htb]
\begin{center}
\begin{tabular}{|c|c|c|c|}
\hline
$n$ & size of sample & $\approx d(n)$ & $\approx d(n)/\sqrt{n}$ \\
\hline
10000 & 2000 & 1.8262177 & 0.018262176 \\
15000 & 2000 & 1.9742892 & 0.016120004 \\
20000 & 3000 & 2.0621564 & 0.014581648 \\
25000 & 4000 & 2.1949573 & 0.013882129 \\
30000 & 5000 & 2.203268 & 0.012720575 \\
35000 & 6000 & 2.3289392 & 0.012448704 \\
40000 & 7000 & 2.3589768 & 0.011794884 \\
\hline
\end{tabular}
\end{center}
\caption{Standard deviation of main diagonal segment}
\label{deviation-center}
\end{table}

\clearpage

\subsection{Average shape in three dimensions}
\begin{figure}[htb]\begin{center}
\includegraphics[scale=0.9]{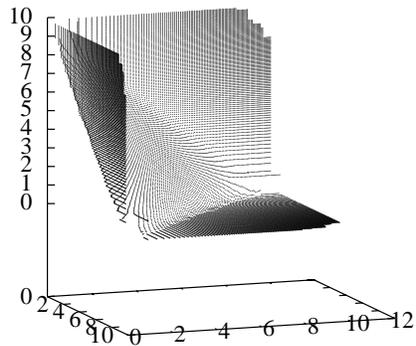}
\caption{Average shape of 400 diagrams of size 10000, drawn in
$(\sqrt{x}, \sqrt{y}, \sqrt{z})$ coordinates}
\label{richardson-shape-3d}
\end{center}\end{figure}

The definition of random Young diagram by Richardson measure can
be easily generalized to three-dimensional case.
There are no known results about the limit shape in this case,
but the figure~\ref{richardson-shape-3d}, obtained by our computations,
leads to assumption that the limit shape complies to
similar equation  $\sqrt{x}+\sqrt{y}+\sqrt{z}=h_3$.

This result shows that the limit shape of 3-dimensional Young
diagrams generated by Richardson process are probably different from
the limit shape for uniformly-distributed diagrams, which was
studied in~\cite{vershik-yakubovich}, and finally
found in~\cite{kenyon-cerf,okounkov}.

\end{document}